\documentclass[12pt,a4paper]{article}
\usepackage{amssymb, amstext, amsmath}
\usepackage{mathtools}
\usepackage{enumerate}

\def\A{\mathcal{A}}
\def\B{\mathcal{B}}
\def\C{\mathcal{C}}
\def\D{\mathcal{D}}
\def\E{\mathcal{E}}
\def\F{\mathcal{F}}
\def\I{\mathcal{I}}

\def\G{\mathcal{G}}

\newtheorem{theorem}{Theorem}[section]
\newtheorem{lemma}{Lemma}[section]

\newtheorem{claim}{Claim}

\title{On the Hilton--Spencer intersection theorems for unions of cycles}
\date{}
\author{Peter Borg\\[5mm]
Department of Mathematics, Faculty of Science, \\ 
University of Malta, Malta\\
\texttt{peter.borg@um.edu.mt}\\
\\
Carl Feghali\\[5mm]
Department of Informatics, University of Bergen, \\ Bergen, Norway\\
\texttt{carl.feghali@uib.no}}

\begin{document}
\maketitle   

\begin{abstract}
A family $\A$ of sets is said to be \emph{intersecting} if every two sets in $\A$ intersect. An intersecting family is said to be \emph{trivial} it its sets have a common element. A graph $G$ is said to be \emph{$r$-EKR} if at least one of the largest intersecting families of independent $r$-element sets of $G$ is trivial. Let $\alpha(G)$ and $\omega(G)$ denote the independence number and the clique number of $G$, respectively. Hilton and Spencer recently showed that if $G$ is the vertex-disjoint union of a cycle ${_*C}$ raised to the power $k^*$ and $s$ cycles ${_1C}, \dots, {_sC}$ raised to the powers $k^*, k_1, \dots, k_s$, respectively, $1 \leq r \leq \alpha(G)$, and
$$\min\big(\omega(_1C^{k_1}), \dots, \omega(_sC^{k_s})\big) \geq 2k^* + 1,$$
then $G$ is $r$-EKR. They had shown that the same holds if ${_*C}$ is replaced by a path and the condition on the clique numbers is relaxed to 
$$\min\big(\omega(_1C^{k_1}), \dots, \omega(_sC^{k_s})\big) \geq k^* + 1.$$
We use the classical Shadow Intersection Theorem of Katona to obtain a short proof of each result for the case where the inequality for the minimum clique number is strict. 
\end{abstract}

\section{Introduction}

Unless stated otherwise, we shall use small letters such as $x$ to denote non-negative integers or elements of a set,
capital letters such as $X$ to denote sets, and calligraphic
letters such as $\mathcal{F}$ to denote \emph{families}
(sets whose members are sets themselves). 
The set of positive integers is denoted by $\mathbb{N}$. The set $\{i \in \mathbb{N} \colon m \leq i \leq n\}$ is denoted by $[m,n]$, $[1,n]$ is abbreviated to $[n]$, and $[0]$ is taken to be the empty set $\emptyset$. For a set $X$, the \emph{power set of $X$} (that is, $\{A \colon A \subseteq X\}$) is denoted by $2^X$. The family of $r$-element subsets of $X$ is denoted by $X \choose r$. The family of $r$-element sets in a family $\mathcal{F}$ is denoted by $\mathcal{F}^{(r)}$. If $\mathcal{F} \subseteq 2^X$ and $x \in X$, then the family $\{A \in \mathcal{F} \colon x \in A\}$ is denoted by $\mathcal{F}(x)$ and called a \emph{star of $\mathcal{F}$} with \emph{centre} $x$. 

A family $\mathcal{A}$ is said to be \emph{intersecting} if for every $A, B \in \mathcal{A}$, $A$ and $B$ intersect (that is, $A \cap B \neq \emptyset$). The stars of a family $\mathcal{F}$ are the simplest intersecting subfamilies of $\mathcal{F}$. We say that $\mathcal{F}$ has the \emph{star property} if at least one of the largest intersecting subfamilies of $\mathcal{F}$ is a star of $\mathcal{F}$.

Determining the size of a largest intersecting subfamily of a given family $\mathcal{F}$ is one of the most popular endeavours in extremal set theory. This started in \cite{EKR}, which features the classical result known as the Erd\H os--Ko--Rado (EKR) Theorem. The EKR Theorem states that if $r \leq n/2$ and $\mathcal{A}$ is an intersecting subfamily of ${[n] \choose r}$, then $|\mathcal{A}| \leq {n-1 \choose r-1}$. 
%
%
Thus, ${[n] \choose r}$ has the star property for $r \leq n/2$ (clearly, for $n/2 < r \leq n$, ${[n] \choose r}$ itself is intersecting). There are various proofs of the EKR Theorem (see \cite{D,FF2,HK2,K,Kat}), two of which are particularly short and beautiful: Katona's \cite{K}, which introduced the elegant cycle method, and Daykin's \cite{D}, using the fundamental Kruskal--Katona Theorem \cite{Ka,Kr}. The EKR Theorem gave rise to some of the highlights in extremal set theory \cite{AK1,F_t1,Kat,W} and inspired many variants and generalizations; 
see \cite{Borg7,DF,F2,F,FT,HST,HT}.

A \emph{graph} $G$ is a pair $(V(G),E(G))$, where $V(G)$ is a set, called the \emph{vertex set of $G$}, and $E(G)$ is a subfamily of ${V(G) \choose 2}$ and is called the \emph{edge set of $G$}. A member of $V(G)$ is called a \emph{vertex of $G$}, and a member of $E(G)$ is called an \emph{edge of $G$}. We may represent an edge $\{v,w\}$ by $vw$. We say that $v$ is \emph{adjacent} to $w$ (in $G$) if $vw$ is an edge of $G$. A subset $I$ of $V(G)$ is an \emph{independent set of $G$} if $vw \notin E(G)$ for every $v, w \in I$. Let $\mathcal{I}_G$ denote the family of independent sets of $G$. An independent set $J$ of $G$ is \emph{maximal} if $J \nsubseteq I$ for each independent set $I$ of $G$ such that $I \neq J$. The size of a smallest maximal independent set of $G$ is denoted by $\mu(G)$. The size of a largest independent set of $G$ is denoted by $\alpha(G)$. A subset $X$ of $V(G)$ is a \emph{clique of $G$} if $vw \in E(G)$ for every $v, w \in X$. The size of a largest clique of $G$ is called the \emph{clique number of $G$} and denoted by $\omega(G)$.

Holroyd and Talbot introduced the problem of determining whether ${\mathcal{I}_G}^{(r)}$ has the star property for a given graph $G$ and an integer $r \geq 1$. Following their terminology, a graph $G$ is said to be \emph{$r$-EKR} if ${\mathcal{I}_G}^{(r)}$ has the star property. The Holroyd--Talbot (HT) Conjecture \cite[Conjecture~7]{HT} claims that $G$ is $r$-EKR if $\mu(G) \geq 2r$. This was verified by Borg \cite{Borg} for $\mu(G)$ sufficiently large depending on $r$ (see also \cite[Lemma~4.4 and Theorem~1.4]{Borgmaxprod}). By the EKR Theorem, the conjecture is true if $G$ has no edges. The HT Conjecture has been verified for several classes of graphs \cite{Borg,Borg1,BH1,BH,FJT,HHS,HS1,HS2,HST,HT,HK,T,Woodroofe}. As demonstrated in \cite{BH}, for $r > \mu(G)/2$, whether $G$ is $r$-EKR or not depends on $G$ and $r$ (both cases are possible). Naturally, graphs $G$ of particular interest are those that are $r$-EKR for all $r \leq \alpha(G)$.

For $n \geq 1$, the graphs $([n], {[n] \choose 2})$ and $([n], \{\{i,i+1\} \colon i \in [n-1]\})$ are denoted by $K_n$ and $P_n$, respectively. For $n \geq 3$, $([n], E(P_n) \cup \{n,1\})$ is denoted by $C_n$. A copy of $K_n$ is called a \emph{complete graph}. A copy $P$ of $P_n$ is called an \emph{$n$-path} or simply a \emph{path}, and a vertex of $P$ is called an \emph{end-vertex} if it is not adjacent to more than one vertex. 
A copy of $C_n$ is called an $n$-cycle or simply a \emph{cycle}. If $H$ is a subgraph of a graph $G$ ($V(H) \subseteq V(G)$ and $E(H) \subseteq E(G)$), then we say that \emph{$G$ contains $H$}. 
For $v, w \in V(G)$, the \emph{distance} $d_G(v,w)$ is $\min\{k \colon v, w \in E(P)$ for some $k$-path $P$ contained by $G\}$. The $k^{\rm th}$ \emph{power of $G$}, denoted by $G^k$, is the graph with vertex set $V(G)$ and edge set $\{vw \colon v,w \in V(G), \, 1 \leq d_G(v,w) \leq k\}$; $G^k$ is also referred to as \emph{$G$ raised to the power $k$}. Note that ${P_n}^k=K_n$ for $k \geq n-1$, and ${C_n}^k=K_n$ for $k \geq n/2$.

The following remarkable analogue of the EKR theorem was obtained by Talbot \cite{T}. 

\begin{theorem}[\cite{T}] \label{thm:talbot}
For $1 \leq r \leq \alpha({C_n}^k)$, ${C_n}^k$ is $r$-EKR. 
\end{theorem}
Talbot introduced a compression technique to prove Theorem~\ref{thm:talbot}. In vague terms, his compression technique rotates anticlockwise the elements of the independent sets of the intersecting family which are distinct from a specified vertex (see Section \ref{section:proof}). 

If $G, G_1, \dots, G_k$ are graphs such that the vertex sets of $G_1, \dots, G_k$ are pairwise disjoint and $G = \left(\bigcup_{i=1}^k V(G_i), \bigcup_{i=1}^k E(G_i) \right)$, then $G$ is said to be the \emph{disjoint union of $G_1, \dots, G_k$}, and $G_1, \dots, G_k$ are said to be \emph{vertex-disjoint}.

Inspired by the work of Talbot, Hilton and Spencer \cite{HS1} went on to prove the following. 

\begin{theorem}[\cite{HS1}] \label{thm:h1}
If $G$ is the disjoint union of one path $P$ raised to the power $k^*$ and $s$ cycles ${_1C}, \dots, {_sC}$ raised to the powers $k_1, \dots, k_s$, respectively, $1 \leq r \leq \alpha(G)$, and 
\begin{equation}\label{condition1}\min\big(\omega(_1C^{k_1}), \dots, \omega(_sC^{k_s})\big) \geq \omega(P^{k^*}), \end{equation}
then $G$ is $r$-EKR. Moreover, for any end-vertex $x$ of $P$, ${\mathcal{I}_G}^{(r)}(x)$ is a largest intersecting subfamily of ${\mathcal{I}_G}^{(r)}$. 
\end{theorem}
The ultimate aim, however, was to obtain a generalization of Theorem~\ref{thm:talbot}, and this was eventually achieved by Hilton and Spencer \cite{HS2} in the following theorem. 

\begin{theorem}[\cite{HS2}] \label{thm:h2}
If $G$ is the disjoint union of $s+1$ cycles ${_*C}, {_1C}, \dots, {_sC}$ raised to the powers $k^*, k_1, \dots, k_s$, respectively, $1 \leq r \leq \alpha(G)$, and 
\begin{equation}\label{condition2} \min\big(\omega(_1C^{k_1}), \dots, \omega(_sC^{k_s})\big) \geq 2k^* + 1, 
\end{equation}
then $G$ is $r$-EKR. Moreover, for any $x \in V(_*C)$, ${\mathcal{I}_G}^{(r)}(x)$ is a largest intersecting subfamily of ${\mathcal{I}_G}^{(r)}$.
\end{theorem}

The proof of Theorem~\ref{thm:h2} is also inspired by Talbot's proof of Theorem~\ref{thm:talbot}. In particular, an essential ingredient in the proof of Theorem \ref{thm:h2} is the use of Theorem~\ref{thm:h1} for the special case where $P$ is a complete graph as the base case of an induction argument. 

In this paper, we give a short proof of Theorem~\ref{thm:h1} and of Theorem~\ref{thm:h2}, except for the cases of equality in conditions (\ref{condition1}) and (\ref{condition2}), respectively. In other words, we prove the following two results. 


\begin{theorem}\label{thm:main}
Theorem~\ref{thm:h1} is true if the inequality in (\ref{condition1}) is strict.
\end{theorem}

\begin{theorem}\label{thm:cycle}
Theorem~\ref{thm:h2} is true if the inequality in (\ref{condition2}) is strict.
\end{theorem}


Our argument is based on the Shadow Intersection Theorem of Katona \cite{Kat}, hence demonstrating yet another application of this classical and useful result in extremal set theory.

\section{The new proof}\label{section:proof}

Let $P, {_1C}, \dots, {_sC}$ be as in Theorem~\ref{thm:h1}. Let $p = |V(P)|$ and $c_i = |V({_iC})|$. 
For $1 \leq i \leq s$, we label the vertices of $_iC$ $1^i, 2^i, \dots, {c_i}^i$ (the superscript $i$ is a label and not a power), where $E({_iC}) = \{\{j^i, (j + 1)^i\} \colon j \in [c_i-1]\} \cup \{n^i,1^i\}$. We may assume that $P = P_p$, that is, $V(P) = [p]$ and $E(P) = \{\{i,i+1\} \colon i \in [p-1]\}$. 
Let $H$ be the union of ${_1C}^{k_1}, \dots, {_sC}^{k_s}$, and let $f: V(H) \rightarrow V(H)$ be the bijection given by
$$f(n^i) = 1^i \quad \mbox{and} \quad f(j^i) = (j + 1)^i \quad \mbox{for $1 \leq i \leq s$ and $1 \leq j \leq c_i$.} $$
Let 
$f^1 = f$, and for any integer $t \geq 2$, let $f^t = f \circ f^{t-1}$ and $f^{-t} = f^{-1} \circ f^{-(t-1)}$. Note that for $t \geq 1$, one can think of $f^t$ as $t$ clockwise rotations, and of $f^{-t}$ as $t$ anticlockwise rotations. For $I \in \mathcal{I}_H$, we denote the set $\{f^t(x) \colon x \in I\}$ by $f^t(I)$, and for $\mathcal{A} \subseteq \mathcal{I}_H$, we denote the family $\{f^t(A) \colon A \in \mathcal{A}\}$ by $f^t(\mathcal{A})$. The notation $f^{-t}(I)$ and $f^{-t}(\mathcal{A})$ is defined similarly.

The new argument presented in this paper lies entirely in the proof of the following important case, which both Theorem~\ref{thm:main} and Theorem~\ref{thm:cycle} pivot on.

\begin{lemma}\label{lem:main}
Theorem~\ref{thm:h1} is true if $P^{k^*}$ is a complete graph and the inequality in (\ref{condition1}) is strict. 
\end{lemma}
Then, Theorem~\ref{thm:main} follows immediately by applying the compression method in \cite{HST}, and Theorem~\ref{thm:cycle} follows by applying the same compression method of Talbot in \cite{T}. 

We now start working towards the proof of Lemma~\ref{lem:main}.

Let $\A$ be a family of $r$-element sets. The shadow of $\A$, denoted by $\partial \A$, is the family $\partial \A = \bigcup_{A \in \A} {A \choose r-1}$. A special case of Katona's Shadow Intersection Theorem \cite{Kat} is that  
\begin{equation}\label{katona} |\A| \leq |\partial \A| \quad \mbox{if $\A$ is intersecting.}
\end{equation}
\textbf{Proof of Lemma \ref{lem:main}.} Suppose that $P^{k^*}$ is a complete graph. Then, $\omega(P^{k^*}) = p$. Suppose $\min\big(\omega(_1C^{k_1}), \dots, \omega(_sC^{k_s})\big) > p$. 
%
%
Note that this implies that for every $i \in [s]$, $c_i \geq p+1$ and, for every $h, j \in [c_i]$ with $h \leq j$,
\begin{equation}\label{condition1''} \mbox{if } j \in \{f^{-q}(h) \colon 1 \leq q \leq p\} \cup \{f^q(h) \colon 1 \leq q \leq p\},\mbox{ then } h^ij^i \in E(_iC^{k_i}).
\end{equation}
%
Let $\A$ be an intersecting subfamily of ${\mathcal{I}_G}^{(r)}$. 
%
%
Recall that $V(P^{k^*}) = [p]$. Let $\A_0 = \{A \in \A \colon A \cap [p] = \emptyset\}$ and $\A_i = \{A \in \A \colon A \cap [p] = \{i\}\}$ for $1 \leq i \leq p$. Since $P^{k^*}$ is a complete graph, $\A_0, \A_1, \dots, \A_p$ partition $\A$. Let $\A_0' = \A_0$ and $\A'_i = \{A \backslash \{i\} \colon A \in \A_i\}$ for $1 \leq i \leq p$. Since $\mathcal{A}$ is intersecting,
\begin{equation} \mbox{for every $i, j \in \{0\} \cup [p]$ with $i \neq j$, each set in $\A_i'$ intersects each set in $\A_j'$.} \label{ref1}
\end{equation}

\begin{claim}\label{claim}
The families $f(\partial\A_0), \A'_1, f^2(\A'_2), \dots, f^p(\A'_p)$ are pairwise disjoint. 
\end{claim}
\textbf{Proof.} 
Suppose $B \in f^i(\A'_i) \cap f^j(\A'_j)$ for some $i, j \in [2,p]$ with $i < j$. Then, $B = f^i(A_i) = f^j(A_j)$ for some $A_i \in \A_i'$ and $A_j \in \A_j'$. Thus, $A_i = f^{j-i}(A_j)$. Since $1 \leq j-i < p$, (\ref{condition1''}) gives us $A_i \cap A_j = \emptyset$, but this contradicts (\ref{ref1}). Similarly, if we assume that $B \in \A'_1 \cap f^i(\A'_i)$ for some $i \in [2,p]$, then we obtain $B \cap A_i = \emptyset$ for some $A_i \in \A_i'$, which again contradicts (\ref{ref1}). Therefore, $\A'_1, f^2(\A'_2), \dots, f^p(\A'_p)$ are pairwise disjoint.

Suppose $B \in f(\partial\A_0) \cap f^i(\A'_i)$ for some $i \in [2,p]$. Then, $B = f(A_0) = f^i(A_i)$ for some $A_0 \in \partial\A_0$ and $A_i \in \A_i'$. Thus, $f^{i-1}(A_i) = A_0 = C \backslash \{x\}$ for some $x \in C \in \A_0$. Since $1 \leq i-1 < p$, (\ref{condition1''}) gives us $A_i \cap C = \emptyset$, but this contradicts (\ref{ref1}). Similarly, if we assume that $B \in f(\partial\A_0) \cap \A'_1$, then we obtain $B \cap C = \emptyset$ for some $C \in \A_0$, which again contradicts (\ref{ref1}). The claim follows.~\hfill{$\Box$} 
\\

Let $\A^*_0 = \{A \cup \{1\} \colon  A \in f(\partial\A_0)\}$, $\A_1^* = \A_1$, and $\A^*_i = \{A \cup \{1\} \colon A \in f^i(\A'_i) \}$ for $2 \leq i \leq p$. For $0 \leq i \leq p$, $\A^*_i \subseteq  {\I_G}^{(r)}(1)$. By Claim~\ref{claim}, $\sum_{i = 0}^p |\A^*_i| \leq |{\I_G}^{(r)}(1)|$. By (\ref{katona}),
$|\A_0| \leq |\partial(\A_0)| = |\A^*_0|$. We have 
$$|\A| = \sum_{i = 0}^p |\A_i| = |\A_0| + \sum_{i = 1}^p |\A_i^*| \leq \sum_{i = 0}^p |\A^*_i| \leq |{\I_G}^{(r)}(1)|,$$
and the lemma is proved.~\hfill{$\Box$}
\\

The full Theorem \ref{thm:main} is now obtained by the line of argument laid out in~\cite{HST}, hence making use of established facts regarding compressions on independent sets.

For any edge $uv$ of a graph $G$, let $\delta_{u,v} \colon \mathcal{I}_G \rightarrow \mathcal{I}_G$ be defined by 
\[\delta_{u,v}(A) = \left\{\begin{array}{ll}
(A \backslash \{v\}) \cup \{u\} & \mbox{if $v \in A, u \notin A$, and $(A \backslash \{v\}) \cup \{u\} \in \mathcal{I}_G$;}\\
A & \mbox{otherwise,} \end{array} \right. \]
and let $\Delta_{u,v} \colon 2^{\mathcal{I}_G} \rightarrow 2^{\mathcal{I}_G}$ be the compression operation defined by
\[\Delta_{u,v}(\mathcal{A}) = \{\delta_{u,v}(A) \colon A \in
\mathcal{A}\} \cup \{A \in \mathcal{A} \colon \delta_{u,v}(A) \in
\mathcal{A}\}.\]
It is well-known, and easy to see, that 
\begin{equation} |\Delta_{u,v}(\mathcal{A})| = |\A|. \nonumber 
\end{equation}
For any $x \in V(G)$, let $N_G(x)$ denote the set $\{y \in V(G) \colon xy \in E(G)\}$. The following is given by \cite[Lemma~2.1]{BH} (which is actually stated for ${\mathcal{I}_G}^{(r)}$ but proved for $\mathcal{I}_G$) and essentially originated in \cite{HST}.

\begin{lemma}[\cite{BH,HST}]\label{BHlemma} If $G$ is a graph, $uv \in E(G)$, $\mathcal{A}$ is an intersecting subfamily of $\mathcal{I}_G$, $\mathcal{B} = \Delta_{u,v}(\mathcal{A})$, $\mathcal{B}_0 = \{B \in \mathcal{B} \colon v \notin B\}$, $\mathcal{B}_1 = \{B \in \mathcal{B} \colon v \in B\}$, and $\mathcal{B}_1' = \{B \backslash \{v\} \colon B \in \mathcal{B}_1\}$, then: \\
(i) $\mathcal{B}_0$ is intersecting; \\
(ii) if $|N_G(u) \backslash (\{v\} \cup N_G(v))| \leq 1$, then $\mathcal{B}_1'$ is intersecting; \\
(iii) if $N_G(u) \backslash (\{v\} \cup N_G(v)) = \emptyset$, then $\mathcal{B}_0 \cup \mathcal{B}_1'$ is intersecting.
\end{lemma}

For a vertex $v$ of a graph $G$, let $G - v$ denote the graph obtained by deleting $v$ (that is, $G-v = (V(G) \backslash \{v\}, \{xy \in E(G) \colon x, y \notin \{v\})$), and let $G \downarrow v$ be the graph obtained by deleting $v$ and the vertices adjacent to $v$ (that is, $G \downarrow v = (V(G) \backslash (\{v\} \cup N_G(v)), \{xy \in E(G) \colon x, y \notin \{v\} \cup N_G(v)\})$).
\\ 
\\
\textbf{Proof of Theorem~\ref{thm:main}.} We use induction on $p$. If $P^{k^*}$ is a complete graph, then the result is given by Lemma~\ref{lem:main}. Note that this captures the base case $p = 1$. Now suppose that $P^{k^*}$ is not a complete graph. Then, $p \geq k^* + 2$. 
If $r=1$, then the result is trivial. Suppose $r > 1$. Let $\mathcal{A}$ be an intersecting subfamily of ${\mathcal{I}_G}^{(r)}$. Let $u = p-1$ and $v = p$. Let $\mathcal{B} = \Delta_{u,v}(\mathcal{A})$, $\mathcal{B}_0 = \{B \in \mathcal{B} \colon v \notin B\}$, $\mathcal{B}_1 = \{B \in \mathcal{B} \colon v \in B\}$, and $\mathcal{B}_1' = \{B \backslash \{v\} \colon B \in \mathcal{B}_1\}$. By Lemma~\ref{BHlemma} (i), $\mathcal{B}_0$ is intersecting. We have $N_G(u) \backslash (\{v\} \cup N_G(v)) = \{p-k^*-1\}$, so, by Lemma~\ref{BHlemma} (ii), $\mathcal{B}_1'$ is intersecting. 
Let $H_0 = P_{p-1}$ and $H_1 = P_{p-k^*-1}$. Clearly, $\mathcal{B}_0 \subseteq {\mathcal{I}_{G-v}}^{(r)}$, $\mathcal{B}_1' \subseteq {\mathcal{I}_{G \downarrow v}}^{(r)}$, $G-v$ is the union of ${H_0}^{k^*}$ and ${_1C}^{k_1}, \dots, {_sC}^{k_s}$, and $G \downarrow v$ is the union of ${H_1}^{k^*}$ and ${_1C}^{k_1}, \dots, {_sC}^{k_s}$. The condition $\min\big(\omega(_1C^{k_1}), \dots, \omega(_sC^{k_s})\big) > \omega(P^{k^*})$ in the theorem gives us $\min\big(\omega(_1C^{k_1}), \dots, \omega(_sC^{k_s})\big) > \omega({H_0}^{k^*}) \geq \omega({H_1}^{k^*})$. By the induction hypothesis, $|\mathcal{B}_0| \leq |{\mathcal{I}_{G-v}}^{(r)}(1)|$ and $|\mathcal{B}_1'| \leq |{\mathcal{I}_{G \downarrow v}}^{(r-1)}(1)|$. We have
\begin{align} |\mathcal{A}| &= |\mathcal{B}| = |\mathcal{B}_0| + |\mathcal{B}_1| \leq |{\mathcal{I}_{G - v}}^{(r)}(1)| + |{\mathcal{I}_{G \downarrow v}}^{(r-1)}(1)| \nonumber \\
&= |\{A \in {\mathcal{I}_G}^{(r)} \colon 1 \in A, v \notin A\}| + |\{A \in {\mathcal{I}_G}^{(r)} \colon 1, v \in A\}| = |{\mathcal{I}_G}^{(r)}(1)|, \nonumber
\end{align}
as required.~\hfill{$\Box$}
\\

We now work towards the proof of Theorem \ref{thm:cycle}. 

Let $G$ consist of the disjoint union of one cycle $_*C$ raised to the power $k^*$ and $s$ cycles ${_1C}, \dots, {_sC}$ raised to the powers $k_1, \dots, k_s$, respectively. Let $c^* = |V(_*C)|$, $c_i = |V(_iC)|$ for $1 \leq i \leq s$, and $n = c^* + \sum_{i = 1}^s c_i = |V(G)|$. As before, for $1 \leq i \leq s$, we label the vertices of $_iC$ $1^i, 2^i, \dots, {c_i}^i$, where $E({_iC}) = \{\{j^i, (j + 1)^i\} \colon j \in [c_i-1]\} \cup \{n^i,1^i\}$.  We shall assume that $_*C = C_{c^*}$, that is, $V(_*C) = [c^*]$ and $E(_*C) = \{\{i,i+1\} \colon i \in [c^*-1]\} \cup \{c^*, 1\}$.
\\
\\
\textbf{Proof of Theorem \ref{thm:cycle}.} We use induction on $c^*$. If $_*C^{k^*}$ is a complete graph, then the result is given by Lemma~\ref{lem:main}. Note that this captures the base case $c^* = 1$. Now suppose that $_*C^{k^*}$ is not a complete graph. Then, $c^* \geq 2k^* + 2$. 
If $r=1$, then the result is trivial. Suppose $r > 1$. Let $\mathcal{A}$ be an intersecting subfamily of ${\mathcal{I}_G}^{(r)}$.

Let $f: V(G) \rightarrow V(G)$ be the Talbot compression \cite{HS2, T} given by
\begin{gather}
f(v) = v \quad \mbox{for $v \in V(G) \backslash V(_*C)$}, \nonumber \\
f(1) = 1, \quad \mbox{and} \nonumber\\
f(1 + j) = 1 + j - 1 \quad \mbox{for $1 \leq j \leq c^* - 1$}. \nonumber
\end{gather}
For $X \in \mathcal{I}_G$ and $\mathcal{X} \subseteq \mathcal{I}_G$, we use the notation $f^t(X)$ and $f^t(\mathcal{X})$ similarly to the way it is used above. Let $F$ be the union of ${C_{c^*-1}}^{k^*}$ and ${_1C}^{k_1}, \dots, {_sC}^{k_s}$. Let $H$ be the union of  ${C_{c^*-k^*-1}}^{k^*}$ and ${_1C}^{k_1}, \dots, {_sC}^{k_s}$. Let 
%
\begin{gather*}
\B = \{A \in \A: 1 \not\in A, \, f(A) \in {\I_F}^{(r)}\},  \\
\C = \{A \in \A: 1 \in A, \, f(A) \in {\I_F}^{(r)}\}, \\
\D_0 = \{A \in \A: 1, \, k^* + 2 \in A\}, \\
\D_i = \{A \in \A: 1 + c^* - i, \, k^* + 2 - i \in A\} \quad \mbox{for $1 \leq i \leq k^*$.}
\end{gather*}
Note that these families partition $\A$. Let
$$
\F = (f^{k^* - 1}(\E) - \{1\}) \cup \bigcup_{i = 0}^{k^*}(f^{k^*}(\D_i) - \{1\}),
$$
where $\E = f(\B) \cap f(\C)$ and, for any family $\G$, $\G - \{1\} = \{G \backslash \{1\}: G \in \G\}$.  

\begin{claim}[see \cite{HS2, T}]\label{claim:final}
The following hold:
\\
(i) $|\A| = |f(\B \cup \C)| + |\F|$; \\
(ii) $f(\B \cup \C)$ is an intersecting subfamily of ${\I_F}^{(r)}$; \\
(iii) $\F$ is an intersecting subfamily of ${\I_H}^{(r-1)}$; \\
(iv) $|{\I_G}^{(r)}(1)| = |{\I_F}^{(r)}(1)| + |{\I_H}^{(r-1)}(1)|$. 
\end{claim}
By the induction hypothesis and Claim \ref{claim:final} (ii)--(iii),  
$|f(\B \cup \C)| \leq |{\I_F}^{(r)}(1)|$ and 
$|\F| \leq |{\I_H}^{(r-1)}|$. Thus, by Claim \ref{claim:final} (i) and Claim \ref{claim:final} (iv), we have
$$
|\A| = |f(\B \cup \C)| + |\F| 
\leq |{\I_F}^{(r)}(1)| + |{\I_H}^{(r-1)}(1)|
=|{\I_G}^{(r)}(1)|,
$$
and the theorem is proved.~\hfill{$\Box$}


\end{document}